\documentclass[11pt]{article}

\usepackage{amsmath,amssymb,enumerate,pb-diagram,framed}
\usepackage[cmtip,arrow]{xy}

\usepackage[applemac]{inputenc}

\def \({\left(}
\def \){\right)}
\def \1{{\bf{1}}}

\def \al{\alpha}

\def \bs{\backslash}

\def \C{{\mathbb C}}

\def \CD{{\cal D}}

\def \CM{{\cal M}}

\def \CO{{\cal O}}
\def \CS{{\cal S}}

\def \df{\ \begin{array}{c} _{\rm def}\\
        ^{\displaystyle =}\end{array}\ }
\def \diag{{\rm diag}}

\def \ds{\displaystyle}
\def\e{\emph}

\def \eqn{\begin{eqnarray*}}
\def \endeqn{\end{eqnarray*}}

\def\f{\fbox}

\def \ga{\gamma}
\def \Ga{\Gamma}

\def \H{{\mathbb H}}
\def \Hom{{\rm Hom}}

\def \Im{{\rm Im}}

\def \la{\lambda}

\def \La{\Lambda}
\def\loc{{\rm loc}}

\def \mod{{\rm mod}}

\def \N{\mathbb N}

\def\ol{\overline}

\def \pa{{\rm par}}
\def \ph{\varphi}
\def \prf{\vspace{5pt}\noindent{\bf Proof: }}

\def \PSL{{\rm PSL}}

\def \qed{\ifmmode\eqno $\square$\else\noproof\vskip
            12pt plus 3pt minus 9pt \fi}
 \def\noproof{{\unskip\nobreak\hfill\penalty50\hskip2em\hbox{}%
     \nobreak\hfill $\square$\parfillskip=0pt%
     \finalhyphendemerits=0\par}}

\def \R{{\mathbb R}}
\def \Re{{\rm Re \hspace{1pt}}}

\def\sm{\smallsetminus}

\def \SL{{\rm SL}}

\def \Z{\mathbb Z}
\def \={\ =\ }

\newtheorem{theorem}{Theorem}[section]

\newtheorem{lemma}[theorem]{Lemma}

\newtheorem{proposition}[theorem]{Proposition}

\renewcommand{\sp}[1]{\left\langle #1 \right\rangle}

\newcommand{\smat}[4]{\(\begin{smallmatrix}#1 & #2 \\ #3 & #4\end{smallmatrix}\)}
\newcommand{\mat}[4]{\(\begin{matrix}#1 & #2 \\ #3 & #4\end{matrix}\)}

\begin{document}
\pagestyle{myheadings} \markright{THE LEWIS CORRESPONDENCE...}

\title{Lewis-Zagier Correspondence for higher order forms}
\author{Anton Deitmar\\ \ \\
Pacific Journal of Mathematics 249.1, 11-21 (2011)}
\date{}
\maketitle

\noindent
{\bf Abstract:}
The Lewis-Zagier correspondence,  which attaches period functions to Maa\ss\ wave forms, is extended to wave forms of higher order, which are higher invariants of the Fuchsian group in question.
The key ingredient is an identification of higher order invariants with ordinary invariants of unipotent twists.
This makes it possible to apply standard methods of automorphic forms to higher order forms.

\tableofcontents

\section*{Introduction}
The Lewis-Zagier correspondence \cite{Lewis,Lewis-Zagier1,Lewis-Zagier2}, see also \cite{Brug}, is a bijection between the space of Maaß wave forms of a fixed Laplace-eigenvalue $\la$ and the space of real-analytic functions on the line satisfying a functional equation which involves the eigenvalue.
The latter functions are called period functions.
In \cite{DeH} this correspondence has been extended to subgroups $\Ga$ of finite index in the full modular group $\Ga(1)$.
One can assume $\Ga$ to be normal in $\Ga(1)$.
The central idea of the latter paper is to consider the action of the finite group $\Ga(1)/\Ga$, and in this way to consider Maaß forms for $\Ga$ as vector-valued Maaß forms for $\Ga(1)$.
This technique can be applied to higher order forms 
\cite{CDO,D-ES,DD, DSr, DKMO,DO,DS} as well, turning the somewhat unfamiliar notion of a higher order invariant into the notion of a classical invariant of a twist by unipotent representation.
This viewpoint has, in the case of Eisenstein series, already been used in \cite{JorgOSull}.
The general framework of higher order invariants and unipotent twists is described in the first section of the present paper.
This way of viewing higher order forms has the advantage that it allows techniques of classical automorphic forms to be applied in the context of higher order forms.
The example of the trace formula will be subject of further investigations by the current author in the near future.
In the present, we apply this technique to extend the Lewis-Zagier correspondence to higher order forms.

We define the corresponding spaces of automorphic forms of higher order in the second section. 
Holomorphic forms of higher order have been defined by various authors.
Maaß forms are more subtle, as it is not immediately clear, how to establish the $L^2$-structure on higher order invariants.
In the paper \cite{DD}, the authors resorted to the obvious $L^2$-structure for the quotient spaces of consecutive higher order forms, which however is unsatisfactory because one whishes to view $L^2$-higher order forms as higher order invariants themselves.
In the present paper this flaw is remedied, as we give a space of locally square-integrable functions on the universal cover of the Borel-Serre compactification whose higher order invariants give the sought for $L^2$-invariants.

In this paper we use a common definition of higher order forms which insists on full invariance under parabolic elements.
In the language of our present Section \ref{Sec1} that means that we take the subgroup $P$ to be the subgroup generated by all parabolic elements.
It is an open question, whether the contents of this paper can be extended to the case $P=\{1\}$, i.e., full higher order invariants.
In this case, Fourier expansions have to be replaced with Fourier-Taylor expansions and thus it is unclear, how the correspondence should be defined.

\section{Higher order invariants and unipotent representations}\label{Sec1}
In this section we describe higher order forms by means of invariants in unipotent representations.

Let $\Ga$ be a group and let $W$ be a $\C[\Ga]$-module.
We here take the field $\C$ of complex numbers as base ring. Most of the general theory works over any ring, but our applications are over $\C$.
Let \fbox{$I_\Ga$} be the augmentation ideal in $\C[\Ga]$ this is the kernel of the \e{augmentation homomorphism} $A:\C[\Ga]\to\C$; $\sum_\ga c_\ga \ga\mapsto\sum_\ga c_\ga$.
The ideal $I_\Ga$ is a vector space with basis $(\ga-1)_{\ga\in\Ga\sm\{ 1\}}$.

In the sequel, we will need  two simple properties of the augmentation ideal which, for the convenience of the reader,  we will prove in the following lemma.
A set $S$ of generators of the group $\Ga$ is called \e{symmetric}, if $s\in S\Rightarrow s^{-1}\in S$.

It is easy to see that 
\begin{itemize}
\item $\C[\Ga]=\C\oplus I_\Ga$.
\item For any given set of generators $S$ of $\Ga$ one has
$$
I_\Ga\=\sum_{s\in S}\C[\Ga](s-1).
$$
\end{itemize}

We also fix a normal subgroup $P$ of $\Ga$.
We let $I_P$ denote the augmentation ideal of $P$ and $\tilde I_P=\C[\Ga] I_P$.
As $P$ is normal, $\tilde I_P$ is a two-sided ideal of $\C[\Ga]$.
For any integer $q\ge 0$ we set
$$
J_q\=I_\Ga^q+\tilde I_P.
$$
The set of $\Ga$-invariants $W^\Ga=H^0(\Ga,W)$ in $W$ can be described as the set of all $w\in W$ with $I_\Ga w=0$.
For $q=1,2,\dots$ we define the set of invariants of type $P$ and order $q$ to be
$$
H^0_{q,P}(\Ga,W)\=\{ w\in W: J_{q}w=0\}.
$$
Then $H_{q,P}^0=H_{q,P}^0(\Ga,W)$ is a submodule of $W$ and we have a natural filtration
$$
0\subset H_{1,P}^0\subset H^0_{2,P}\subset\dots\subset H_{q,P}^0\subset\dots,
$$
and as $I_\Ga H_{q,P}^0\subset H_{q-1,P}^0$, the group $\Ga$ acts trivially on $H_{q,P}^0/H_{q-1,P}^0$.

A representation $(\eta,V_\eta)$ of $\Ga$ on a complex vector space $V_\eta$ is called a \e{unipotent length $q$ representation}, if $V_\eta$ has a $\Ga$-stable filtration
$$
0\subset V_{\eta,1}\subset\dots\subset V_{\eta,q}=V_\eta,
$$
such that $\Ga$ acts trivially on each quotient $V_{\eta,k}/V_{\eta,k-1}$ where $k=1,\dots, q$ and $V_{\eta,0}=0$.

Let a unipotent length $q$ representation $(\eta,V_\eta)$ be given.
We also assume that it is $P$-trivial, i.e., the restriction to the subgroup $P$ is the trivial representation.
There is a natural map 
$$
\Phi_\eta:\Hom_\Ga(V_\eta,W)\otimes V_\eta\to W
$$
given by $\al\otimes v\mapsto\al(v)$.

\begin{lemma}
Let $W$ be a $\C[\Ga]$-module.
The submodule $H_{q,P}^0(\Ga,W)$ constitutes a $P$-trivial, unipotent length $q$ representation of $\Ga$.
If the group $\Ga$ is finitely generated, then the space $H_{q,P}^0(\Ga,W)$ is the sum of all images $\Phi_\eta$ when $\eta$ runs over the set of all $P$-trivial, unipotent length $q$ representations which are finite dimensional over $\C$.
\end{lemma}

\prf
The first assertion is clear.
Assume now that $\Ga$ is finitely generated.
The space $H_{q,P}^0(\Ga,W)$ needn't be finite dimensional.
We use induction on $q$ to show that for each $w\in H_{q,\pa}^0(\Ga,W)$ the complex vector space $\C[\Ga]v$ is finite dimensional.
For $q=1$ we have $\C[\Ga]w=\C w$ and the claim follows.
Next let $w\in H_{q+1,P}^0(\Ga,W)$ and let $S$ be a finite set of generators of $\Ga$.
Then
$$
\C[\Ga]w\=\C w+I_\Ga w\=\C w+\sum_{s\in S}\C[\Ga](s-1)w.
$$
The element $(s-1)w$ lies in $H_{q,P}^0(\Ga,W)$, so by induction hypothesis, the claim follows.
\qed

Assume from now on that $\Ga$ is finitely generated.
The philosophy pursued in the rest of the paper is this:
\begin{framed}
\begin{center}
\e{Once you know $\Hom_\Ga(V_\eta,W)$ for every $P$-trivial,\\ finite dimensional unipotent length $q$ representation,\\ you know the space $H_{q,P}^0(\Ga,W)$.}
\end{center}
\end{framed}
So, instead of investigating $H^0_{q,P}(\Ga,W)$ one should rather look at
$$
\Hom_\Ga(V_\eta,W)\cong (V_\eta^*\otimes W)^\Ga,
$$
which is often easier to handle.
In fact, it is enough to restrict to a generic set of $\eta$.
As an example of this philosophy consider the case $q=2$.
For each group homomorphism $\chi:\Ga/P\to (\C,+)$ one gets a $P$-trivial, unipotent length $q$ representation \f{$\eta_\chi$} on $\C^2$ given by
$\eta_\chi(\ga)\=\smat 1{\chi(\ga)}\ 1$.

We introduce the following notation
$$
\bar H_{q,P}^0\=\bar H_{q,P}^0(\Ga,W)\=  H_{q,P}^0(\Ga,W)/ H_{q-1,P}^0(\Ga,W)\=   H_{q,P}^0/ H_{q-1,P}^0.
$$

\begin{proposition}
The space $ H_{2,P}^0(\Ga,W)$ is the sum over all images $\Phi_{\eta_\chi}$, where $\chi$ ranges in $\Hom(\Ga/P,\C)\sm \{0\}$.
For any two $\chi\ne\chi'$ one has
$$
\Im(\Phi_{\eta_\chi})\cap \Im(\Phi_{\eta_{\chi'}})\=  H^0(\Ga,W).
$$
In other words, one has
$$
\bar H_1^0\= \bigoplus_{\chi}\Im(\Phi_{\eta_\chi})/ H^0.
$$
\end{proposition}

\prf
We make use of the \e{order lowering operator}
$$
\La: \bar H_{q,P}^0\to\Hom(\Ga/P,\bar H_{q-1,P}^0)\cong\Hom(\Ga/P,\C)\otimes\bar  H_{q-1,P}^0,
$$
where the last isomorphism is due to the fact that $\Ga$ is finitely generated.
This operator is defined as
$$
\La(w)(\ga)\= (\ga -1)w.
$$
One sees that this indeed is a homomorphism in $\ga$ by using the fact that 
$$
(\ga\tau-1)\equiv (\ga-1)+(\tau-1)\ \mod\ I^2
$$ 
for any two $\ga,\tau\in\Ga$.
The map $\La$ is clearly injective.

Let now $w\in \Im(\Phi_{\eta_\chi})\cap \Im(\Phi_{\eta_{\chi'}})$ for $\chi\ne\chi'$.
Then
$$
\La(w)\in\chi\otimes H^0\cap\chi'\otimes H^0,
$$
and the latter space is zero as $\chi\ne\chi'$.
For surjectivity, let $w\in H_1^0$.
Then $\La(w)=\sum_{i=1}^n\chi_i\otimes w_i$ with $w_i\in H^0$, and so $w\in\sum_{i=1}^n\Im(\phi_{\eta_{\chi_i}})$.
\qed

\section{Higher order forms}\label{Sec2}
We next define some spaces of automorphic forms of higher order, like holomorphic modular forms or Maaß wave forms \cite{CDO,D-ES,DD, DSr, DKMO,DO,DS}.
The holomorphic case has been treated in various other places.
Maaß forms are more subtle, as it is not immediately clear, how to establish the $L^2$-structure on higher order invariants.
In the paper \cite{DD} the authors resorted to the obvious $L^2$-structure for the quotient spaces $\bar H_{q,P}^0$, which however is unsatisfactory because one whishes to view $L^2$-higher order forms as higher order invariants themselves.
In the present paper this flaw is remedied, as we give a space of locally square-integrable functions on the universal cover of the Borel-Serre compactification whose higher order invariants give the sought for $L^2$-invariants.
We also give a guide how to set up higher order $L^2$-invariants in more general cases, like general lattices in locally compact groups, when there is no such gadget as the Borel-Serre compactification around.
In that case, Lemma \ref{Lem2.1} tells you how to define the $L^2$-structure once you have chosen a fundamental domain for the group action.

Let \f{$G$} denote the group $\PSL_2(\R)=\SL_2(\R)/\{\pm 1\}$.
It has the group \f{$K$}$={\rm PSO}(2)={\rm SO}(2)/\{\pm 1\}$ as a maximal compact subgroup.
Let
\f{$\Ga(1)$}$=\PSL_2(\Z)$ be the full modular group.
Let \f{$\Ga$}$\subset\Ga(1)$ be a normal subgroup of finite index which is torsion-free.
For every cusp $c$ of $\Ga$ fix \f{$\sigma_c$}$\in\Ga(1)$ such that $\sigma_c\infty= c$ and
$\sigma_c^{-1}\Ga_c\sigma_c=\pm\smat 1{N_c\Z}\ 1$.
The number \f{$N_c$}$\in\N$ is uniquely determined and is called the \e{width} of the cusp $c$.
Let \f{$\H$}$=\{ z\in\C: \Im(z)>0\}$ be the upper half plane and let \f{$\CO(\H)$} be the set of holomorphic functions on $\H$.
We fix a weight $k\in 2\Z$ and define a (right-) action of $G$ on functions $f$ on $\H$  by
$$
f|_k\ga(z)\=
(cz+d)^{-k}f\(\frac{az+b}{cz+d}\),\qquad \ga=\mat abcd.
$$
We define the space \f{$\CO_{\Ga,k}^M(\H)$} to be the set of all $f\in\CO(\H)$, such that for every cusp $c$ of $\Ga$ the function $f|_k\sigma_c$ is, in the domain $\{\Im(z)>1\}$, bounded by a constant times $\Im(z)^A$ for some $A>0$.

Further we consider the space \f{$\CO_{\Ga,k}^S(\H)$} of all $f\in\CO(\H)$, such that for every cusp $c$ of $\Ga$ the function $f|_k\sigma_c$ is, in the domain $\{\Im(z)>1\}$, bounded by a constant times $e^{-A\Im(z)}$ for some $A>0$.

These two spaces are preserved not only by $\Ga$, but also by the action of the full modular group $\Ga(1)$.

The normal subgroup $P$ of $\Ga$ will be the subgroup \f{$\Ga_\pa$} generated by all parabolic elements.
We then write $H_{q,\pa}^0$ for $H_{q,P}^0$.
We consider the space of modular functions of weight k and order $q$,
$$
\text{\f{$M_{k,q}(\Ga)$}}\= H_{q,\pa}^0(\Ga,\CO_{\Ga,k}^M(\H)),
$$
as well as the corresponding space of cusp forms
$$
\text{\f{$S_{k,q}(\Ga)$}}\= H_{q,\pa}^0(\Ga,\CO_{\Ga,k}^S(\H)).
$$
Then every $f\in M_{k,q}(\Ga)$ possesses a \e{Fourier expansion}  at every cusp $c$ of the form
$$
f|_k\sigma_c(z)\=\sum_{n=0}^\infty a_{c,n}\, e^{2\pi i\frac n{N_c}z}.
$$
A function $f\in M_{k,q}(\Ga)$ belongs to the subset $S_{k,q}(\Ga)$ if and only if $a_{c,0}=0$  for every cusp $c$ of $\Ga$.

As the group $\Ga$ is normal in $\Ga(1)$, the latter group acts on the finite dimensional spaces $M_{k,q}(\Ga)$ and $S_{k,q}(\Ga)$.
These therefore give examples of finite dimensional representations of $\Ga(1)$ which become unipotent length $q$ when restricted to $\Ga$.

By a \emph{Maa\ss\ wave form} for the group $\Ga$  and parameter $\nu\in\C$ we mean a
function $u\in L^2(\Ga\bs\H)$ which is twice continuously differentiable and
satisfies
\begin{eqnarray*}
 \Delta u&=&\({\textstyle\frac 14 -\nu^2}\) u.
\end{eqnarray*}
By the regularity of solutions of elliptic differential equations this condition implies that $u$ is real analytic.
Let \f{$\CM_\nu$}$=\CM_\nu(\Ga)$ be the space of all Maa\ss\ wave forms for $\Ga$.
Note that sometimes in the definition of Maa\ss\ forms, instead of the $L^2$-condition, a weaker condition on the growth at the cusps is imposed.

Next we define Maa\ss-wave forms of higher order.
First we need the higher order version of the Hilbert space $L^2(\Ga\bs \H)$.
For this recall the construction of  the Borel-Serre compactification $\ol{\Ga\bs\H}$ of $\Ga\bs\H$, see \cite{BS}.
First one constructs a space \f{$\H_\Ga$}$\supset\H$ by attaching to each cusp $c$ of $\Ga$ a real line to $\H$ and then one equips this set with a suitable topology such that $\Ga$ acts properly discontinuously and the quotient $\Ga\bs\H_\Ga$ is the Borel-Serre compactification.
The space $\H_\Ga$ is constructed in such a way, that for a given (closed) fundamental domain  $D\subset\H$ of $\Ga\bs\H$ which has finitely many geodesic sides, the closure $\ol{D}$ in $\H_\Ga$ is a fundamental domain for $\Ga\bs\H_\Ga$.
By the discontinuity of the group action, this has the following consequence: For every compact set $K\subset\H_\Ga$ there exists a finite set $F\subset \Ga$ such that $K\subset F\ol D=\bigcup_{\ga\in F}\ga\ol D$.

Now we extend the hyperbolic measure to $\H_\Ga$ in such a way that the boundary $\partial \H_\Ga=\H_\Ga\sm\H$ is a nullset.
Let \f{$L^2_\loc(\H_\Ga)$} be the space of local $L^2$-functions on $\H_\Ga$.
Then $\Ga$ acts on $L^2_\loc(\H_\Ga)$.
Since $\Ga$ acts discontinuously with compact quotient on $\H_\Ga$, one has 
$$
L_\loc^2(\H_\Ga)^\Ga\= L^2(\Ga\bs\H).
$$
Define
 the space \f{$L^2_q(\Ga\bs\H)$} as the space of all $f\in L^2_\loc(\H_\Ga)$ such that $J_{q}f=0$, so in other words,
$$
L^2_q(\Ga\bs\H)\=H_{q,\pa}^0\(\Ga,L^2_\loc(\H_\Ga)\).
$$
Then $L_1^2(\Ga\bs\H)=L^2(\Ga\bs\H)$ is a Hilbert space in a natural way.
We want to install Hilbert space structures on the spaces $L^2_q(\Ga\bs\H)$ for $q\ge 2$ as well.
For this purpose we introduce the space \f{$F_q$} of all measurable functions $f:\H\to\C$ such that $J_{q}f=0$ modulo nullfunctions.
Then $L_q^2(\Ga\bs\H)$ is a subset of $F_q$.

\begin{lemma}\label{Lem2.1}
Let $S\subset\Ga$ be a finite set of generators which is assumed to be symmetric and to contain the unit element.
Let $D\subset\H$ be a closed fundamental domain of $\Ga$ which has finitely many geodesic sides.

Any $f\in F_q$ is uniquely determined by its restriction to
$$
S^{q-1}D\=\bigcup_{s_1,\dots,s_{q-1}\in S}s_1\dots s_{q-1} D.
$$
One has 
$$
F_q\cap L^2(S^{q-1}D)\= L^2_q(\Ga\bs\H),
$$
where on both sides we mean the restriction to $S^qD$, which is unambiguous by the first assertion.
In this way the space $L^2_q(\Ga\bs\H)$ is a closed subspace of the Hilbert space $L^2(S^{q-1}D)$.
The induced Hilbert space topology on $L^2_q(\Ga\bs\H)$ is independent of the choices of $S$ and $D$, although the inner product is not.
The action of the group $\Ga(1)$ on $L_q^2(\Ga\bs\H)$ is continuous, but not unitary unless $q=1$.
\end{lemma}

\prf
We have to show that any $f\in L_q=L_q^2(\Ga\bs\H)$ which vanishes on $S^{q-1}D$, is zero.
We use induction on $q$.
The case $q=1$ is clear.
Let $q\ge 2$ and write \f{$\bar L_q$}$= L_q/L_{q-1}$.
Consider the \e{order lowering operator}
$$
\La:L_q\to\Hom(\Ga,\bar L_{q-1})\cong\Hom(\Ga,\C)\otimes \bar L_{q-1},
$$
given by
$$
\La(f)(\ga)\= (\ga -1)f.
$$
The kernel of $\La$ is $L_{q-1}$.
Now assume $f(S^{q-1}D)=0$.
Then for every $s\in S$ we have $(s-1)f(S^{q-2}D)=0$ and hence, by induction hypothesis, we conclude $(s-1)f=0$.
But as $S$ generates $\Ga$ this means that $\La(f)=0$ and so $f\in L_{q-1}$, so, again by induction hypothesis, we get $f=0$.

We next show that
$$
F_q\cap L^2(S^{q-1}D)\=F_q\cap L^2(S^{q-1+j}D)
$$ for every $j\ge 0$.
The inclusion ``$\supset$'' is clear.
We show the other inclusion by induction on $q$ and $j$.
For $q=1$ or $j=0$ there is no problem.
So assume the claim proven for $q$.
Let $f\in F_{q+1}\cap L^2(S^{q+j}D)$ and let $s\in S$.
Then $f(sz)=f(z)+f(sz)-f(z)$, the function $f(z)$ in in $L^2(S^{q+j}D)$ and the function $f(sz)-f(z)$ is in $F_{q}\cap L^2(S^{q+j-1}D)= L^2(S^{q+j}D)$ by induction hypothesis.
It follows that $f\in L^2(s^{-1}S^{q+j}D)$ and since this holds for every $s$ we get $f\in L^2(S^{q+j+1}D)$ as claimed.

We now come to
$$
F_q\cap L^2(S^{q-1}D)\= L^2_q(\Ga\bs\H).
$$
Let $f\in F_q\cap L^2(S^{q-1}D)$. For every compact subset $K$ of $\H_\Ga$ there exists $j\ge 0$ such that $K\subset S^{q-1+j}D$.
Therefore, $f$ is in $L^2(K)$ for every compact subset $K$ of $H_\Ga$.
As the latter space is locally compact, $f$ is in $L^2_\loc(\H_\Ga)$.
Since $I_\Ga^{q}f=0$ we get $f\in L_q^2(\Ga\bs\H)$.
For the other inclusion let $f\in L_q^2(\Ga\bs\H)$.
As $S^{q-1}D$ is relatively compact in $\H_\Ga$ it follows that $f\in L^2(S^{q-1}D)$ as claimed.

We next show independence of the topology of $S$.
So let $S'$ be another set of generators.
Then there exists $l\in\N$ such that $S'\subset S^l$.
Hence ist suffices to show that the topology from the inclusion $L_q^2(\Ga\bs\H)\subset L^2(S^{q-1}D)$ coincides with the topology from the inclusion $L_q^2(\Ga\bs\H)\subset L^2(S^{q+j}D)$ for every $j\ge 0$.
If a sequence tends to zero in the latter, it clearly tends to zero in the first.
The other way round is proven by induction on $j$ similar to the above.
In particular, the continuity of the $\Ga(1)$-action follows.

Finally, we show the independence of $D$.
Let $D'$ be another closed fundamental domain with finitely many geodesic sides.
Then there exists $l\in\N$ such that $D'\subset S^lD$ and the claim follows along the same lines as above.
\qed

We define the space \f{$\CM_{\nu,q}$}$=\CM_{\nu,q}(\Ga)$ of \e{Maa\ss-wave forms of order $q$} to be the space of all $u\in L^2_q(\Ga\bs\H)$ which are  twice continuously differentiable and satisfy
$$
\Delta u\= \(\frac14-\nu^2\)u.
$$
Fix a finite dimensional representation $(\eta,V_\eta)$ of $\Ga(1)$, which is $\Ga_\pa$-trivial and becomes a unipotent length $q$ representation on restriction to $\Ga$.

We  set \f{$\CM_{\nu,q,\eta}$} equal to $\( V_\eta\otimes\CM_{\nu,q}\)^{\Ga(1)}$. 
Likewise we define \f{$\tilde\CM_{\nu,q}(\Ga)$}$=\tilde\CM_{\nu,q}$
as the space of all $u\in F_q(\Ga)$ which 
are twice continuously differentiable and 
satisfy $\Delta u\= \(\frac14-\nu^2\)u$, and 
we set \f{$\tilde\CM_{\nu,q,\eta}$}$=(V_\eta\otimes\tilde\CM_{\nu,q})^{\Ga(1)}$.

\begin{lemma}
Let \f{$\CD_\nu'$} be the space of all distributions $u$ on $\H$ with $\Delta u=(\frac14-\nu^2)u$.
Then
$$
\tilde\CM_{\nu,q,\eta}
\= (V_\eta\otimes\tilde\CM_{\nu,q})^{\Ga(1)}
\= (V_\eta\otimes \CD_\nu')^{\Ga(1)}
$$
and
$$
\CM_{\nu,q,\eta}
\= (V_\eta\otimes\CM_{\nu,q})^{\Ga(1)}
\=\(V_\eta\otimes (\CD_\nu'\cap L^2_\loc(\H_\Ga))\)^{\Ga(1)}.
$$
\end{lemma}

\prf
The inclusion ``$\subset$'' is obvious in both cases.
We show ``$\supset$''.
In the first case, the space on the left can be described as the space of all smooth functions $u:\H\to V_\eta$ satisfying $\Delta u=(\frac14-\nu^2)u$ as well as $J_{q+1}u=0$ and $u(\ga z)=\eta(\ga)u(z)$ for every $\ga\in\Ga(1)$.
Now let $u\in (V_\eta\otimes \CD_\nu')^{\Ga(1)}$.
As $u$ satisfies an elliptic differential equation with smooth coefficients, $u$ is a smooth function with $\Delta u=(\frac14-\nu^2)u$.
The condition $u(\ga z)=\eta(\ga)u(z)$ is clear.
Finally, the condition $J_{q+1}u=0$ follows from that, as $\eta|_\Ga$, being $\Ga_\pa$-trivial and  unipotent of length $q$, satisfies $\eta(J_{q+1})=0$.
Hence the first claim is proven.
The second is similar.
\qed

As in the holomorphic case, every Maaß-form $f\in\CM_{\nu,q}(\Ga)$ has a Fourier expansion at every cusp $c$,
$$
f(\sigma_c z)\=\sum_{n=0}^\infty a_{c,n}(y) \, e^{2\pi i \frac{n}{N_c} x},
$$
with smooth functions $a_{c,n}(y)$.

We define the space \f{$\CS_{\nu,q}$} of \emph{Maa\ss\ cusp forms} to be the
space of all 
$f\in\CM_{\nu,q}$ with $a_{c,0}(y)=0$ for every cusp $c$.
We also set \f{$\CS_{\nu,q,\eta}$}$=(V_\eta\otimes\CS_{\nu,q})^{\Ga(1)}$.

Note that since $\eta$ is unipotent of length $q$ on $\Ga$, we have
$$
\CS_{\nu,q,\eta}\=\(V_\eta\otimes\CS_{\nu,q'}\)^{\Ga(1)}
$$
for every $q'\ge q$.

\section{Setting up the transform}\label{Sec3}

It is the aim of this note to extend the Lewis Correspondence  
\cite{DeH,Lewis, Lewis-Zagier1, Lewis-Zagier2} to the case of higher order forms.
We will explain the approach in the case of cusp forms first.

\begin{center}
Throughout, let $(\eta,V_\eta)$ be a finite dimensional representation of $\Ga(1)$ which becomes $\Ga_\pa$-trivial and unipotent\\ of length $q$ when restricted to $\Ga$.
\end{center}

We fix the following notation for the canonical
generators of $\Ga(1)$:
$$
S\=\pm\mat 01{-1}0,\quad {\rm and} \quad T\=\pm\mat 1101.
$$
Then $S^2=\1=(ST)^3$, and $T$ is of infinite order.
Note that there exists $N\in\N$ such that $T^N\in\Ga$ as $\Ga(1)/\Ga$ is a finite group.
Let $N$ be minimal with this property, then $N=N_\infty$ is the width of the cusp $\infty$ of $\Ga$.
We then have $\eta(T)^N=\eta(T^N)=1$, as $\eta$ is trivial on parabolic elements of $\Ga$.

Let \f{$\Psi_{\nu,\eta}$} be the space of all holomorphic functions $\psi\colon
\C\sm(-\infty,0]\to V_\eta$ satisfying the \e{Lewis equation}
\begin{equation}\label{LewisEquation}
\eta(T)\psi(z)=
\psi(z+1)+(z+1)^{-2\nu-1}\eta(ST^{-1})\psi
\left(\frac z{z+1}\right)
\end{equation}
and the \e{asymptotic formula}
\begin{equation}\label{LimitEquation}
\begin{array}{l}
0\ =\   \ds e^{{+}\pi i \nu}\lim_{\Im(z)\to\infty}
             \left(\psi(z)+z^{-2\nu-1}\eta(S)\psi\left(\frac {-1}z\right)\right)+\\
\phantom{0\ =\ } \ds  +e^{{-}\pi i\nu}\lim_{\Im(z)\to-\infty}
             \left(\psi(z)+z^{-2\nu-1}\eta(S)\psi\left(\frac{-1}z\right)\right),\\
   \end{array}
\end{equation}
where both limits are supposed to exist.

Let {$A$} denote the subgroup of $G$
consisting of diagonal matrices and let {$N$} be the
subgroup of upper triangular matrices with $\pm 1$ on the diagonal.
The group $G$ then as a manifold is a direct product
$G=ANK$. For $\nu\in\C$ and
$a=\pm\diag(t,t^{-1})\in A$, $t>0$, let $a^\nu = t^{2\nu}$.
We insert the factor $2$ for compatibility reasons.

Let \f{$(\pi_\nu,V_{\pi_\nu})$} denote the
principal series representation of $G$ with parameter $\nu$. The
representation space $V_{\pi_\nu}$ is the Hilbert space of all functions
$\ph\colon G\to \C$ with
$\ph(anx)=a^{\nu+\frac 12}\ph(x)$ for $a\in A, n\in N, x\in G$,
and $\int_K |\ph(k)|^2\, dk <\infty$
modulo nullfunctions. The representation is
$\pi_\nu(x)\ph(y)=\ph(yx)$.
There is a special vector \f{$\ph_0$} in $V_{\pi_\nu}$ given by
$$
\ph_0(ank)\= a^{\nu+\frac 12}.
$$
This vector is called the \emph{basic spherical function} with parameter
$\nu$.

For a continuous $G$-representation $(\pi,V_\pi)$ on a topological vector space $V_\pi$ let \f{$\pi^\omega$}
denote the subrepresentation on the space of analytic vectors, i.e. $V_{\pi^\omega}$ consists of all
vectors $v$ in $V_\pi$ such that for every continuous linear map $\al \colon V_\pi\to\C$ the map
$g\mapsto\al(\pi(g)v)$ is real analytic on $G$.
This space comes with a natural topology. Let \f{$\pi^{-\omega}$} be its topological dual.
In the case of $\pi=\pi_\nu$ it is known that $\pi_\nu^\omega$ and $\pi_\nu^{-\omega}$ are in perfect duality,
i.e., they are each other's topological duals.
The vectors in $\pi_\nu^{-\omega}$ are called \emph{hyperfunction vectors} of the representation $\pi_\nu$.

As a crucial tool we will use the space
\begin{equation*}
\text{\f{$A_{\nu,\eta}^{-\omega}$}}=(\pi_\nu^{-\omega}\otimes\eta)^{\Ga(1)}=
H^0(\Ga(1),\pi_\nu^{-\omega}\otimes\eta)
\end{equation*}
and call it the space of \emph{$\eta$-automorphic hyperfunctions}.

For an automorphic hyperfunction $\al\in A_{\nu,\eta}^{-\omega}$ we consider the
function $u\colon G\to V_\eta$ given by
$$
u(g) \ \df\  \sp{\pi_{-\nu}(g)\ph_0,\al}.
$$
Here $\sp{\cdot,\cdot}$ is the canonical pairing $\pi_{-\nu}^\omega\times \pi_{-\nu}^{-\omega}\otimes\eta\to V_\eta$.
Then $u$ is right $K$-invariant, hence can be viewed as a function on $\H$.
As such it lies in $\tilde\CM_{\nu,\eta}$ since $\al$ is $\Ga$-equivariant and
the Casimir operator on $G$, which induces $\Delta$, is scalar on $\pi_\nu$
with eigenvalue $\frac 14-\nu^2$.
The transform \f{$P$}$\colon \al\mapsto u$ is called the \emph{Poisson transform}. It
follows from \cite{Schlichtkrull}, Theorem 5.4.3, that the Poisson transform
\begin{equation*}
P\colon A_{\nu,\eta}^{-\omega}\ \to\ \tilde\CM_{\nu,\eta}
\end{equation*}
is an isomorphism for $\nu\not\in \frac 12+\Z$.

For $\al\in A_{\nu,\eta}^{-\omega}$ put
$$
\text{\f{$\psi_\al(z)$}}\ = \
f_\al(z)-z^{-2\nu-1}\eta(S)f_\al\left(\frac{-1}z\right),
$$
with $f_\al$ 
such that the function $z\mapsto
(1+z^2)^{\nu+\frac12} f_\al(z)$
represents the restriction $\al|_\R$.
Then the \emph{Bruggeman transform}
$B\colon\al\mapsto\psi_\al$ maps
$A_{\nu,\eta}^{-\omega}$ to $\Psi_{\nu,\eta}$.
It is a bijection if $\nu\notin \frac 12+\Z$, as can be seen similar to \cite{DeH}, Proposition 2.2.

For $\nu\not\in \frac 12+\Z$ we finally define the \emph{Lewis transform} as the
map \f{$L$}$\colon \CM_{\nu,\eta}\to \Psi_{\nu,\eta}$,
given by
\begin{equation*}
L\ \df\  B\circ P^{-1}.
\end{equation*}

\begin{theorem}\label{Lewistrafobijektiv}
{\rm (Lewis transform; cf. \cite{Lewis-Zagier2}, Thm.~1.1)}
For $\nu\not\in \frac 12+\Z$ and $\Re(\nu)> -\frac 12$
the Lewis transform is a bijective linear map from
the space of Maa\ss\ cusp
forms $\CS_{\nu,q,\eta}$ to the space $\Psi_{\nu,\eta}^o$ of period functions.
\end{theorem}

\prf The proof runs, with small obvious changes, along the lines of the corresponding result \cite{DeH}, Theorem 3.3.
\qed

{\small
\noindent
Universit\"at T\"ubingen,
Mathematisches Institut,  
Auf der Morgenstelle 10, 
72076 T\"ubingen, Germany, 
{\tt deitmar@uni-tuebingen.de}
}

\small
\begin{thebibliography}{XXX}

\bibitem{BS}
\bf Borel, A.; Serre, J.-P.:
\it Corners and arithmetic groups. Avec un appendice: Arrondissement des variétés à coins, par A. Douady et L. Hérault. 
\rm Comment. Math. Helv. 48 (1973), 436--491.

\bibitem{Brug}
\bf Bruggeman, R.W.:
\it Automorphic forms, hyperfunction cohomology, and period functions.
\rm J. reine angew. Math. 492 (1997), 1--39.
\bibitem{CDO} 
\bf Chinta, G.; Diamantis, N.; O'Sullivan, C.: 
\it{Second order modular forms} 
\rm Acta Arith., 103 (2002), 209-223.
\bibitem{DeH}
\bf Deitmar, Anton; Hilgert, Joachim:
\it A Lewis Correspondence for submodular groups. 
\rm Forum Math. 19, no. 6, 1075-1099 (2007).

\bibitem{D-ES}
\bf Deitmar, A.:
\it Higher order group cohomology and the Eichler-Shimura map. 
\rm  J. reine u. angew. Math. 629, 221-235 (2009).

\bibitem{DD}
\bf Deitmar, A.; Diamantis, N.:
\it Automorphic forms of higher order.
\rm Journal of the LMS. doi:10.1112/jlms/jdp015 (2009).

\bibitem{DSr}
\bf Diamantis, N.; Sreekantan, R.:
\it Iterated integrals and higher order automorphic forms.
\rm Commentarii Mathematici Helvetici 81 (2006), 481--494.

\bibitem{DKMO}
\bf Diamantis, N.; Knopp, M.; Mason, G.; O'Sullivan, C.:
\it $L$-functions of second-order cusp forms.
\rm Ramanujan J. 12 (2006), no. 3, 327--347.

\bibitem{DO}
\bf Diamantis, N.; O'Sullivan, C.:
\it The dimensions of spaces of holomorphic second-order automorphic forms and their cohomology.
\rm Trans. Amer. Math. Soc. 360, no. 11, 5629-5666 (2008).

\bibitem{DS}
\bf Diamantis, N.; Sim, D.:
\it The classification of higher-order cusp forms.
\rm J. Reine Angew. Math. 622, 121-153 (2008).

\bibitem{JorgOSull}
\bf Jorgenson, Jay; O'Sullivan, Cormac:
\it Unipotent vector bundles and higher-order non-holomorphic Eisenstein series.
\rm J. Théor. Nombres Bordeaux 20 (2008), no. 1, 131--163.

\bibitem{Lewis}
\bf Lewis, J.:
\it Spaces of holomorphic functions equivalent to the even
Maass cusp forms.
\rm Invent. Math. 127  (1997), 271--306.

\bibitem{Lewis-Zagier1}
 \bf Lewis, J.; Zagier, D.:
 \it Period functions and the Selberg zeta function for the modular group.
 \rm The mathematical beauty of physics (Saclay, 1996), 83--97, Adv. Ser. Math.
Phys., 24, World Sci. Publishing, River Edge, NJ, 1997.

\bibitem{Lewis-Zagier2}
\bf Lewis,
J.; Zagier, D.:
\it Period functions for Maass wave forms.
\rm I. Ann. of Math. (2) 153 (2001), 191--258.

\bibitem{Schlichtkrull}
\bf Schlichtkrull, H.:
\it Hyperfunctions and harmonic analysis on symmetric spaces.
\rm Progress in Mathematics, 49. Birkh\"auser
Boston, Inc., Boston, MA, 1984.


\end{thebibliography}
\end{document}